\documentclass[letterpaper, 10 pt, conference]{ieeeconf}

\IEEEoverridecommandlockouts        

\usepackage{cite}
\usepackage{amsmath,amssymb,amsfonts}
\usepackage{algorithmic}
\usepackage{graphicx}
\usepackage{textcomp}
\usepackage{xcolor}
\usepackage{verbatim}
\usepackage{ esint }
\usepackage{ dsfont }
\usepackage{ amssymb }
\usepackage{graphicx,subcaption}
\usepackage{algorithmic}
\usepackage{multirow}
\usepackage{comment}
\usepackage{graphics} 

\newtheorem{assumption}{Assumption}
\newtheorem{theorem}{Theorem}
\newtheorem{proposition}{Proposition}

\title{\LARGE \bf
Building Temperature Control: A Distributed Escort Dynamical approach
}

\author{M. Sawant, J. Moyalan, J. Koonamparampath, A. Sheikh, S. Wagh, and N. Singh \\ EED, VJTI, Mumbai, India
\thanks{M. Sawant, J. Moyalan, J. Koonamparampath, A. Sheikh, S. Wagh and N. Singh are with the Electrical Engineering Department, Veermata Jijabai Technological Institute, Mumbai 400019, INDIA
        {\tt\small mssawant\_m17@ee.vjti.ac.in}}
}

\begin{document}

\maketitle

\thispagestyle{empty}
\pagestyle{empty}

\begin{abstract}
    The constrained multi-agent optimization problem of distributed resource allocation is addressed using the evolutionary game theoretic framework. The issue of building temperature control is analyzed in which the controller is to devise a scheme to distribute available scarce power to every room to regulate their temperature as per the user's comfort in the best possible manner. The paper correlates the global constraint of fixed resource amount with the constant population size. The respective population game is evaluated by means of a dynamical model of the evolutionary game theory to find the necessary control action. The robustness of optimal solution with respect to minor fluctuations in the temperature distribution is characterized using evolutionary stable strategy (ESS). Along with the global constraint, the problem formulation constitutes local constraint over an individual control unit located in every room. The classical dynamical models of evolutionary game theory such as replicator dynamics, logit dynamics, etc. fail to incorporate respective constraints. With the escort evolutionary dynamical (ED) model it is possible to address these local constraints through the concept of the intersection of simplices. However, evaluation of these classical dynamics along with ED is driven by expected payoff obtained by the overall population, which renders a centralized implementation approach. To mitigate this central dependency a distributed version of the ED model referred to as distributed escort dynamics (DED) is proposed. The control action devised adopting DED approach is shown to provide smooth trajectory tracking along with the low start-up transience when compared with distributed interior point (DIP) method.
\end{abstract}
\begin{table}[ht]
\caption{List of Symbols}
\begin{tabular}{|l|l|}
\hline
\multirow{2}{0.7cm}{$N$}  & Total number of zones within the building \\
 & ( rooms and walls )\\
\hline
$\mathcal{R}$& Set comprising of rooms\\\hline
$\mathcal{W}$& Set comprising of walls\\\hline
$\mathcal{Z}$& Set representing all the zones\\\hline
$\alpha_{i,j}$& Thermal condunctance between $i^{th}$ and $j^{th}$ entity\\\hline
$u_i$& Output power of $i^{th}$ actuator unit\\\hline
$U$& Total power generation at every instance\\\hline
$k$ & Total number of rooms\\
\hline
$\boldsymbol{t^{set}}$ & Set point temperature profiles\\
\hline
$\boldsymbol{t^{a}}$ & Ambient temperature profile\\
\hline
\multirow{2}{0.7cm}{$\boldsymbol{d^i}$} & Disturbance profile causing temperature deviations \\
 & at $i^{th}$ instance\\
\hline
\multirow{2}{0.7cm}{$\boldsymbol{u^{lo}}$} & Lower bounds on the individual actuator outputs \\
 & ( local constraints )\\
\hline
\multirow{2}{0.7cm}{$\boldsymbol{u^{up}}$} & Upper bounds on the individual actuator outputs \\
 & ( local constraints )\\
\hline
\multirow{2}{0.7cm}{$\mathcal{G}$} & Graph ensuring information transfer among \\
 & the neighbouring actuator mechanisms\\
\hline
$\boldsymbol{t^i}$ & Actual temperature profiles of zones at $i^{th}$ instance\\
\hline
$\boldsymbol{x^i}$ & Population distribution at $i^{th}$ instance\\
\hline
$f_j^i$ & Payoff obtained for employing $j^{th}$ strategy at $i^{th}$ instance\\
\hline
$\boldsymbol{x^{lo}}$& Upper bounds on the population proportions \\
\hline
$\boldsymbol{x^{up}}$& Lower bounds on the population proportions\\
\hline
$\boldsymbol{\eta}$ & Escort function accommodating $\boldsymbol{x^{lo}}$\\
\hline
$\boldsymbol{\xi}$ & Escort function accommodating $\boldsymbol{x^{up}}$\\
\hline
\end{tabular}
\end{table}

\section{Introduction}
The rapid speed of urbanization has resulted in an ever-increasing number of commercial as well as residential building infrastructures. The energy efficient operation of these buildings is very crucial not only from financial but also from environmental perspective \cite{nejat2015global},\cite{lechtenbohmer2011potential}. The variety of building maintenance and regulation processes are taken care of by using heating, ventilation and air conditioning (HVAC) systems \cite{yang2014thermal},\cite{perez2008review}. Of these, the process of thermal regulation results in the majority of energy consumption.

In building temperature control (BTC) formulation \cite{obando2015distributed}, the main goal of an HVAC system is to regulate the temperature within the rooms as per respective residents' comfort, while optimizing total energy consumption. The limited availability of the energy presents BTC issue as a multi-objective, multi-agent resource allocation problem. In this, the control mechanism installed within every room acts as a different agent whose objective is to regulate respective room temperature around resident's desired setpoint. Mainly there are two ways to address such multi-agent optimization problems, i.e. centralized and distributed. 

Though centralised approach in which central entity gathers system information and provides an exact allocation scheme, it is a very time-consuming process. Also, the dependency on a central authority raises security and privacy concerns. Moreover, it is not always possible to satisfy the requirement of central authority in distributed scenario \cite{murphey2000approximate}. On the other hand, the distributed approach, instead of depending on a global state, solves the optimization problem based on the locally available information \cite{chevaleyre2006issues} of the respective agent's neighbourhood. The distributed method provides consensus-based approach \cite{obando2015distributed} in which a fair play is established ensuring equal payoff to every agent.

The paper exploits the interacting multi-agent perspective to propose an evolutionary game theory based escort dynamical (ED) approach \cite{harper2011escort} to address BTC problem. Unlike classical evolutionary dynamics, the ED model is shown to accommodate bounds on the individual control dynamics \cite{ovalle2016escort} along with global resource constraint through the intersection of simplices. Moreover, based on a distributed version of the classical evolutionary model discussed in \cite{barreiro2016distributed}, a consensus-based distributed implementation of ED model referred to as distributed escort dynamics (DED) model is proposed. This model with the property of positive invariantness is shown to converge to an equilibrium state corresponding to equal payoff to every agent.

The performance evaluation of DED implementation is carried out through comparative analysis with widely used resource allocation protocol \cite{obando2015distributed}, distributed interior point (DIP) method based BTC problem implementation. It has been shown that DED approach provides smooth trajectory tracking along with a low start-up transience as compared to the DIP approach. The major contributions of the proposed research work are:
\begin{itemize}
    \item Implementation of an evolutionary game theoretical based approach to address resource allocation problem of BTC.
    \item Distributed version ED model has been proposed and convergence analysis of the same has been provided.
\end{itemize}

The rest of the paper is structured as follows. Section II describes the BTC process and provides the control objective. In section III, ED dynamics have been introduced, also, the distributed version of ED, referred to as DED is proposed. Section IV analyzes the BTC problem implementation using DED. Representative case studies and concluding remarks are present thereafter.

\section{Mathematical Prerequisite and Notations}
The vector quantities are represented using bold lowercase symbols. The symbolic representation $x_j^i$, where $i, j$ belongs to class of whole numbers, represents $j^{th}$ component of vector $\boldsymbol{x}$ at $i^{th}$ instance. The identity matrix is represented as $I$. The difference notation $|\cdot|$ represents the Euclidean norm.

\subsection{Graph Theory}\label{AA}
In the multi-agent system considered in this paper, modeling of the communication network with the help of graph allows the agents to coordinate their decisions as given in \cite{godsil2013algebraic}. A graph is described by the triplet $\mathcal{C}=( \mathcal{S},\mathcal{L},\mathcal{A})$. $\mathcal{S}={\{1,. . ., k}\}$ represents the set of nodes, $\mathcal{L}\subseteq\mathcal{S}\times\mathcal{S}$ represents the set of edges connecting the nodes and matrix $\mathcal{A}$ represents a $k\times k$ nonnegative matrix whose elements satisfy: $\rho_{ij} = 1$ if and only if $(i,j)\in\mathcal{L}$, and $\rho_{ij}=0$ if and only if $(i,j)\notin\mathcal{L}$. The nodes and edges of the graph corresponds to agents and communication channels of the multi-agent system respectively. The neighbours of node $i$ i.e. the set of nodes that are able to receive/send information from/to node $i$ is represented by $\mathcal{N}_{i}={\{j \in \mathcal{S} : (i,j) \in \mathcal{L}}\}$.

The graphical models utilized throughout the text are assumed to constitute several properties, i.e.
\begin{assumption}
The graph does not contain self loops, $ \rho_{ii} = 0 \text{ } \forall i\in \mathcal{S}$ i.e.
\end{assumption}
\begin{assumption}
The graph is undirected. $\rho_{ij}=\rho_{ji} \text{ } \forall i, j \in \mathcal{S}$
\end{assumption}
\begin{assumption}
The graph is connected, i.e. for every pair of nodes in $\mathcal{S}$ there exists either direct or indirect path which connects the two.
\end{assumption}

For the graph $\mathcal{C}$, the Laplacian matrix $L(\mathcal{C})$ is a $k \times k$ dimensional matrix with individual entities defined as,
\begin{equation} 
l_{ij} = 
\begin{cases}
\sum_{j \in \mathcal{S}} \rho_{ij}, \quad &\text{if} \quad i = j\\
-\rho_{ij}, \quad \quad &\text{if} \quad  i \neq j\\
\end{cases}
\end{equation}

\section{Building Temperature Control process}
\begin{figure}[ht]
\centering
\includegraphics[width=0.45\textwidth]{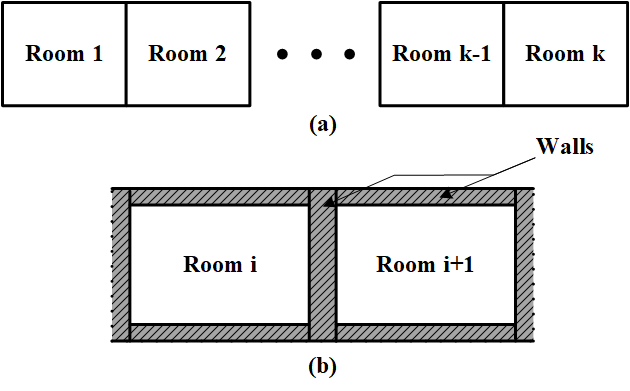}
\caption{(a) Building architecture with sequential room arrangement (b) Top view representation of two adjacent rooms surrounded with walls}
\label{Building Architecture}
\end{figure}

A building can be considered to be partitioned mainly into two sets, a set of $k$ rooms $\mathcal{R}$ and $m$ walls $\mathcal{W}$ also referred to as zones Fig. \ref{Building Architecture}. Let $\mathcal{Z}$ defines the composition of sets $\mathcal{R}$ and $\mathcal{W}$ i.e. the total number of zones. Every room has a desired temperature profile $\boldsymbol{t^{set}}$ defined over a time period by its resident referred to as the setpoint temperature profile. The temperature within every room is controlled by a central thermal unit (CTU) \cite{bravo2015distributed}. In this case, a ventilation system connected to CTU is spread throughout the building and is also connected to every room. The CTU controls the airflow which affects the room temperature through an actuator mechanism installed therein.

\subsection{The Temperature Dynamics}
The building temperature evaluated over a certain time duration exhibits a dynamic behaviour. At any time instance, the temperature within any zone gets directly affected by temperature of adjacent zones. Also, the ambient temperature of the surrounding affects the overall building temperature \cite{nguyen2014relationship} by directly interfering with the temperature of zones which are in direct contact with it. The effect is mainly dependent on the value of thermal conductance between the two interacting bodies. Also, the impromptu opening and closing of windows and doors, and resident's body temperature affects the temperature of the respective zone.

The combined effect of all the factors on the temperature of the zone is represented through a dynamical equation as,
\begin{equation}
    \theta_i \dot{t}_i = \sum_{j=1}^N\alpha_{i, j}(t_j - t_i) + \alpha_{i,a}(t_i^a - t_i) + v_i(u_i + d_i), \forall i\in \mathcal{Z}\label{temp_dynamics}
\end{equation}
where, $\alpha_{i,j}$ is the thermal conductance between $i^{th}$ and $j^{th}$ zone. Its value is positive if respective zones are in direct contact with each other, else it is zero. Similarly, $\alpha_{i,a}$ represents the thermal conductance between outside environment and zone (walls) which is in direct contact with it. $\theta_i$ represents the thermal capacitance value for the $i^{th}$ zone while $N$ is the total number of zones. 

The control parameter $u_i$ corresponds to the output of the actuator installed in $i^{th}$ zone. The parameter $v_i$ is binary variable, with possible values between $\{0, 1\}$. Positive value of $v_i$ ensures that the $i^{th}$ zone that is being evaluated in indeed a room. The thermal disturbance such as impromptu opening and closing of windows and/or effect of body temperature of resident with respect to the $i^{th}$ zone is characterized as $d_i$. The detailed description of BTC process can be found in \cite{bravo2015distributed}.

\subsection{The Control Architecture}
The CTU generates a fixed $U$ amount of power at every instance which is to be distributed among actuator mechanisms so as to modulate respective rooms' temperature around desired setpoints.
\begin{equation}
    \sum_{i=1}^k u_i = U\label{old_global_constraint}
\end{equation}

However, in this case, because of an actuator output constraints, the individual actuator output $u_i$ is bounded within upper and lower bounds.
\begin{align}
    u_i^{lo} \leq u_i \leq u_i^{up} ,\quad \forall i \in \mathcal{R} \label{local_constraint}
\end{align}

Depending on the constraints value \eqref{local_constraint}, the \eqref{old_global_constraint} may or may not hold true. To accommodate this uncertainty, a positive semidefinite variable $u_{k+1}$ is introduced, so that,
\begin{equation}
    \sum_{i=1}^{k+1} u_i = U\label{global_constraint}
\end{equation}
In the BTC framework, the variable $u_{k+1}$ corresponds to the unused power at respective time instance. 

\subsection{BTC Problem Definition}
The building represents a global system which comprises of various subsystems i.e. different zones. These zones interact with each other according to dynamical model \eqref{temp_dynamics}. The global objective here is to utilize optimal power to regulate individual room temperature around its setpoint value. 
\begin{equation}
    \text{minimize } \sum_{i=1}^k\frac{(t_i - t_i^{set})^2}{2},\quad \forall i \in \mathcal{A} \label{global_objective}
\end{equation}

The global objective function \eqref{global_objective} can be decoupled and solved locally at individual subsystem level, i.e.
\begin{equation}
    \text{minimize } \frac{(t_i - t_i^{set})^2}{2},\quad \forall i \in \mathcal{A} \label{local_objective}
\end{equation}
Hence, the task is to develop individual actuator implementation strategies to attain the local objective \eqref{local_objective}, while satisfying the global \eqref{global_constraint} and local constraints \eqref{local_constraint}, respectively.

\section{BTC as a Dynamics Distributed Resource Allocation problem}
Unlike general resource allocation problem formulations where objective function $f_i$ is real valued function of resource variable $x_i$, the payoff formulation in BTC problem is dependent on the state \eqref{local_objective} i.e. temperature values of rooms. These states are driven through a dynamical model \eqref{temp_dynamics} which is affected by the resource being allocated \eqref{plant}. 
\begin{equation}
D^P_i :
\begin{cases}
\dot t_i = g(\boldsymbol{t}, u_i) \\
f_i = h(\boldsymbol{t}, u_i)\\
\end{cases}
\label{plant}
\end{equation}
Here $D^P_i$ represents the state dynamical model which is driven by instantaneous state values $\boldsymbol{t}$ and control parameter which also corresponds to resource $u_i$. 

\begin{figure}[htbp]
\centerline{\includegraphics[scale=0.6]{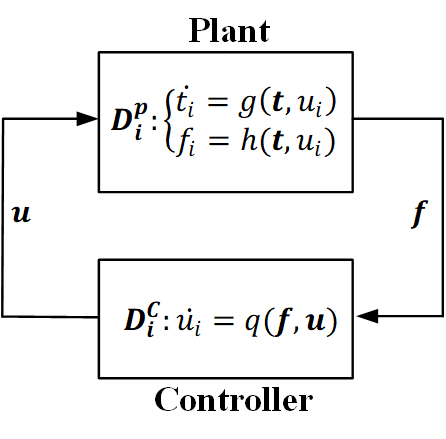}}
\caption{General system schematic}
\label{general_schematic}
\end{figure}

The system model is interconnected with the controller module $D^C_i$ which perform the resource allocation task according to a function $u$, \eqref{controller}.
\begin{equation}
D^C_i : \dot x_i = q(f_i,f_j), \forall j \in \mathcal{N}_i
\label{controller}
\end{equation}
Fig. \ref{general_schematic} provides the schematic representation of overall system. In DED approach to BTC formulation, functions $g(\cdot), h(\cdot), \text{ and } q(\cdot)$ corresponds to \eqref{temp_dynamics}, \eqref{local_objective}, and \eqref{DED}, respectively.

The dynamical system model \eqref{plant} is considered to be reached equilibrium point when output consensus is attained \cite{obando2015distributed}. When a set of subsystems is considered to be reached output consensus if $lim_{t\rightarrow \inf}|f_i - f_j| = 0$, for all $i, j = 0, 1, 2, ..., k$, where $f_i$ is the output of $i^th$ subsystem i.e. room.

Hence, to obtain the steady state, the controller model \eqref{controller} must drive the \eqref{plant} to the output consensus. Also, controller parameter should satisfy system defined constraints such as global resource constraint \eqref{global_constraint} and local bounds \eqref{local_constraint}.

\subsection{Convergence Analysis to the Output Consensus}

Let $(\boldsymbol{t^*}, \boldsymbol{u^*})$ be the rest point of the dynamical model Fig. \ref{general_schematic}. To analyze the system with respect to its convergence to an output consensus, the system dynamics \eqref{plant}, \eqref{controller} are represented in error coordinates with respect to $(\boldsymbol{t^*}, \boldsymbol{u^*})$.
\begin{equation}
ErrD^P_i :
\begin{cases}
\dot e_{t_i} = e_g(\boldsymbol{e_t}, \boldsymbol{e_u}) \\
e_{f_i} = e_h(\boldsymbol{e_t}, \boldsymbol{e_u})\label{error_plant_module}
\end{cases}
\end{equation}
\begin{equation}
ErrD^C_i : \dot{e_{u_i}} = e_q(\boldsymbol{e_f}, \boldsymbol{e_u})\label{error_controller_module}
\end{equation}

The feedback interconnection Fig. \ref{general_schematic} is considered to satisfy the uniqueness assumption. The statement is as follows;

\begin{assumption}
For a feedback interconnection of plant \eqref{plant} and controller \eqref{controller} dynamics represents in error coordinates with respect to $(\boldsymbol{t^*}, \boldsymbol{u^*})$ as \eqref{error_plant_module}, \eqref{error_controller_module}, if $e_g(0, \boldsymbol{e_u}) = 0,$ then $e_{u_i} = 0, \forall i = 1, 2, ..., k.$
\label{unique}
\end{assumption}

The \textit{Assumption 4} ensures the uniqueness of the rest point $(\boldsymbol{t^*}, \boldsymbol{u^*})$ i.e. feedback interconnection dynamics Fig. \ref{general_schematic} will only have one convergence point. For the feedback interconnection of multi-agent systems satisfying the \textit{Assumption 4}, the \textit{Theorem 3.2.1} from \cite{obando2015distributed} provides a set of sufficient conditions which ensures the convergence of feedback interconnection Fig. \ref{general_schematic} to output consensus. The reformulated theorem is stated as,

\begin{theorem}
For a feedback interconnection of system dynamics \eqref{plant}, \eqref{controller} satisfying the \textit{Assumption 4}, let $(\boldsymbol{t^*}, \boldsymbol{u^*})$ be an equilibrium point, where control law $q(\cdot)$ is of the form \eqref{controller}, if following conditions hold true:
\begin{enumerate}
    \item [C1.] The communication graph $\mathcal{G}$ for multi-agent system \eqref{controller} is connected.
    \item [C2.] The plant dynamics \eqref{plant} expressed in error coordinates with respect to $(\boldsymbol{t^*}, \boldsymbol{u^*})$, \eqref{error_plant_module}, is strictly passive from the input $e_u$ to the output $e_f$ with radially unbounded storage function.
    \item [C3.] $u^*$ and $u(0)$ satisfy the global resource constraint \eqref{global_constraint}.
\end{enumerate}
Then, \eqref{plant} reaches output consensus.
\end{theorem}

The condition C1 ensures that there always exists a path between two distinct agents to facilitate information transfer. The passive nature of plant dynamics is ensured by the fulfillment of condition C2. Moreover, the interconnection of two passive systems results in the convergence to an equilibrium point. In order to utilize this property of passive interconnection not only plant dynamics \eqref{plant} but also controller dynamics \eqref{controller} must be passive. According to the proposition provided in \cite{obando2015distributed}, the passive nature of controller dynamics is ensured by the fulfillment of the condition C3. The proposition is reformulated as follows,

\begin{proposition}
Given that the $\boldsymbol{u^*}$ satisfies the constant population requirement \eqref{global_constraint}, if $\boldsymbol{u(0)}$ satisfies \eqref{global_constraint} and graph $\mathcal{G}$ is connected, then multi-agent dynamics in error coordinates \eqref{error_controller_module} is passive and lossless from the input $\boldsymbol{e_f}$ to the output $\boldsymbol{-e_u}$.
\end{proposition}

In this paper, an escort dynamics based distributed algorithm whose interconnection with the BTC dynamics \eqref{temp_dynamics} as shown in Fig. \ref{general_schematic} satisfies the conditions enlisted in \textit{ Theorem 1}, providing convergence to an output consensus.

\section{Escort Evolutionary Game Dynamics}
The evolutionary game theory analyzes the interactions among the large number of players with comparatively less number of strategies. Instead of analyzing individual player dynamics, it evaluates the evolution of the spread of available strategies over the entire population. During this evolution process, every player within the population chooses a certain strategy according to some classical dynamical revision protocols such as, replicator dynamics, logit dynamics, projection dynamics, etc. These models, also referred to as classical dynamical models are of the form. 
\begin{equation}
\dot x_i^{next}= \mathcal{D}(x_i^{present}, f_i^{present}) \label{D}
\end{equation}
These dynamical models drives the population proportions towards best population proportion, such that,
\begin{equation}
    x_i^*f_i^* > x_if_i \text{ for all } i = 1, 2, ..., k\label{ess_property}
\end{equation}
where, $\boldsymbol{x} = [ x_1, x_2, ..., x_k]^T$ represents population state vector. $x_i$ corresponds to the proportion of population employing $i^{th}$ pure strategy. Similarly, $f_i$ corresponds to the  payoff received for employing $i^{th}$ strategy. $\boldsymbol{x^*}$ and $\boldsymbol{f^*}$ represent best possible population proportion vector and payoff vector, respectively. 

\begin{figure}[ht]
\centering
\includegraphics[width=0.35\textwidth]{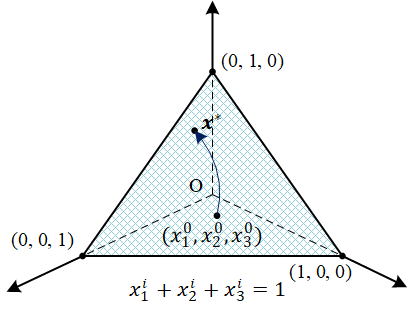}
\caption{Simplex geometry in three dimensional Euclidean space}
\label{Standard Simplex}
\end{figure}

The normalized population proportion vector $\boldsymbol{x}$ can be considered as a population distribution on a $k-1$ dimensional statistical simplex manifold.
\begin{equation}
\Delta_k = \{\boldsymbol x\in \mathds{R}^{k}|\sum_{i = 1}^{k}{x_i} = 1, x_i\geq 0\}\label{simplex}
\end{equation}
hence, the evolution towards best possible population distribution corresponds to the motion of initial population state vector $\boldsymbol{x_0}$ over $(k-1)$ dimensional manifold of $\Delta_k$, \eqref{simplex}, as shown in Fig. \ref{Standard Simplex}.

In the perspective of BTC process, the population proportion vector $\boldsymbol{x}$ corresponds to the actuator output profile $\boldsymbol{u}$. Also, the global resource constraint \eqref{global_constraint} gets reflected in simplex formulation \eqref{simplex}. In evolutionary game theoretical framework, the local bounds on the actuator mechanism \eqref{local_constraint} are correlated as lower and upper bounds on individual strategy proportion \eqref{local_simplex_bound}.
\begin{equation}
\Gamma = \{x_i^{lo}\leq x_i\leq x_i^{up}\}, \text{ for all } i = 1, 2, ..., k\label{local_simplex_bound}
\end{equation}
Then the task is to determine best possible population distribution corresponding to the best output profile for individual actuator mechanism.

However, the classical dynamical models of the form \eqref{D} do not have any provisions to accommodate local individual proportion constraints \eqref{local_simplex_bound}. In \cite{barreiro2018constrained}, the mixture of classical dynamical models is considered to incorporate such proportion bound. In this paper, the evolutionary game theoretical model of ED \cite{harper2011escort} has been used to address this issue. 

\subsection{Escort Dynamics Features}
The ED dictates the evolution of normalised population distribution over $k$ pure strategies under the influence of payoff achieved \cite{harper2011escort}. The continuous time representation of the ED is,
\begin{equation}
    \dot{x_i} = \phi(x_i)(f_i(\boldsymbol{x}) - f_{\phi})\text{ for all } i = 1, 2, ..., k \label{ED}
\end{equation}
where, $\phi(x_i)$ is the positive semidefinite function of population proportion, also referred to as escort function. Also, $f_i(\boldsymbol{x})$ and $f_{\phi}$ represents the payoff obtained for employing $i^{th}$ pure strategy and weighted average payoff, respectively. 

The weighted average payoff $f_{\phi}$ is computed as, 
\begin{equation}
    f_{\phi} = \frac{1}{\Phi(\boldsymbol{x})}\sum_{i=1}^{k}\phi_k(x_i)f_i(\boldsymbol{x})\label{average_payoff}
\end{equation}
where, $\Phi(\boldsymbol{x}) = \sum_{i=1}^k{\phi_i(x_i)f_i(\boldsymbol{x})}$. The presence of $f_{\phi}$ within dynamical model \eqref{ED} renders centralised implementation approach. For positive escort functions, \eqref{average_payoff} can be considered as the expected value of the payoff function over a probability distribution defined by escort functions as,
\begin{equation}
    \hat{\phi}(\boldsymbol{x}) = \frac{1}{\Phi{\boldsymbol{x}}}[\phi(x_1), \phi(x_2), ..., \phi(x_k)]\label{escort_distribution}
\end{equation}
summation of \eqref{ED} over entire population for all $k$ pure strategies is computed as
\begin{equation}
\sum_{i = 1}^{k}{\dot x_i} = \sum_{i = 1}^{k}{x_if_i(\boldsymbol x)} - f_{\phi}\sum_{i = 1}^{k}{x_i} = 0 \label{ED constant population}
\end{equation}
that is initially defined population size, $\sum_{i = 1}^{k}{x_i} = 1$, remains constant throughout the evolution for $t>0$. Hence, the definition of escort function is very crucial to constrain population proportion within more restrictive region.

However, the presence of weighted average term $f_{\phi}$ within ED \eqref{ED} demands that at every revision instance the player must know the payoff received by every other player belonging to population. This center oriented arrangement is quite restrictive in nature. Also, corresponding real life implementation will require more communicational as well as computational infrastructure. To mitigate this central dependency the paper proposed a graph theoretic based distributed version of ED referred to as DED.

\subsection{Distributed Escort Dynamics formulation}
The central dependency of ED \eqref{ED} model is removed through a graph theoretical based approach. In this case, at every revision instance the $i^{th}$ strategy proportion evolution is governed by the local information obtained from its neighbourhood $\mathcal{N}_i$. The proposed DED model is represented as,
\begin{equation}
    \dot{x}_i = \phi(x_i)\sum_{j\in\mathcal{N}_i}{\phi(x_j)[f_j(\boldsymbol{x}) - f_i(\boldsymbol{x})]} \text{ for all } i = 1, 2, ..., k \label{DED}
\end{equation}

Here, $\phi(\cdot)$ is the escort function. Instead of depending on the expected average payoff as in \eqref{ED}, the DED dynamics $i^{th}$ proportion \eqref{DED} is driven by the payoffs obtained by neighbourhood proportion. The DED model can be considered as consensus based algorithm. i.e. the DED dynamics \eqref{DED} reaches steady state either if respective $\phi(\cdot)$ becomes zero or if the payoff values obtained by respective proportion is same as that of its neighbouring proportions.

The DED dynamics \eqref{DED} can also be represented as,
\begin{equation}
    \dot{x}_i = \sum_{j = 1}^{k}\rho_{i,j}[f_j(\boldsymbol{x}) - f_i(\boldsymbol{x})]\text{ for all } i = 1, 2, ..., k \label{DED_version}
\end{equation}
where, $\rho_{i,j}$ governs the connectivity between the two distinct population proportions. Its value represents the weight of the links between two adjacent population proportions. Here, the connectivity graph $\mathcal{G}$ is assumed to be undirected, hence $\rho_{i,j}$ equals $\rho_{j,i}$. The link weightage $\rho_{i,j}$ is formulated as,
\begin{equation}
    \begin{aligned}
        \rho_{i,j} &= \phi(x_i)\phi(x_j) &\text{ if } (i, j) \in \Sigma\\
        \rho_{i, j} &= 0 &\text{ Otherwise }
    \end{aligned}
\end{equation}

\subsubsection{The Positive Invariantness of DED}
Similar to \eqref{ED constant population}, the summation of \eqref{DED_version} over all strategy proportion is evaluated as,
\begin{equation}
\begin{aligned}
    \sum_{i=1}^{k}\dot{x}_i = \sum_{i=1}^k\bigg[\sum_{j = 1}^{k}\rho_{i,j}[f_j(\boldsymbol{x}) - f_i(\boldsymbol{x})]\bigg] =\\
    \rho_{1,1}(f_1 - f_1) + \rho_{1,2}(f_2 - f_1) + ... +\rho_{1,k}(f_k - f_1)\\
        +\rho_{2,1}(f_1 - f_2) + \rho_{2,2}(f_2 - f_2) + ... +\rho_{2,k}(f_k - f_2)\\
        :\\
            +\rho_{k,1}(f_1 - f_k) + \rho_{k,2}(f_2 - f_k) + ... +\rho_{k,k}(f_k - f_k) \label{positive invariantness}
    \end{aligned}
\end{equation}

As the connectivity graph $\mathcal{G}$ is undirected the expression \eqref{positive invariantness} equates to zero. i.e.
\begin{equation}
    \sum_{i=1}^k \dot{x}_i = 0\label{positive_inv_proved}
\end{equation}
equation \eqref{positive_inv_proved} ensures that the population size remains constant throughout the evolution process. i.e. if for the initial distribution $\sum_{i=1}^{k} x_i = 1$ then it will remain the same for all $t > 0$. This property is referred to as positive invariantness.

The property of positive invariantness limits the possible population distribution vectors within the manifold corresponding to $\sum_{i=1}^{k} x_i = 1$. However, to restrict the population dynamics within more constrained bounds \eqref{local_simplex_bound} the escort formulation is very crucial. In order to accommodate system defined local constraints \eqref{local_constraint} on actuator output which correlates into bounded nature of individual proportion coefficient $x_i$ for all $i = 1, 2, ..., k$, \eqref{local_simplex_bound} the escort function formulation is discussed.

\subsection{Escort function formulation for BTC problem}
The geometric representation of individual strategy proportion bounds \eqref{local_simplex_bound} corresponds to the intersection of two $k-1$ dimensional simplices, $\Delta_k^{lo}$ and $\Delta_k^{up}$, where
\begin{align}
\Delta_k^{lo} &= \{\boldsymbol x \in \mathds{R}^k | x_i \geq {x_i^{lo}}, \sum_{i = 1}^{k} x_i = 1\}\label{low simplex}\\
\Delta_k^{up} &= \{\boldsymbol x \in \mathds{R}^k | x_i \leq {x_i^{up}}, \sum_{i = 1}^{k} x_i = 1\}\label{up simplex}
\end{align} 

Structures $\Delta_k^{lo}$ and $\Delta_k^{up}$ along with $\Delta_k$ are convex in nature which allow them to represent just with the help of their vertices. Let $S^{lo}$ and $S^{up}$ be $k \times k$ dimensional matrices whose column vectors represent vertices of $\Delta_k^{lo}$ and $\Delta_k^{up}$, respectively. Construction of $S^{lo}$ and $S^{up}$ is given as,
\begin{align}
S^{lo} &= X^{lo} + \sigma^{lo}I\label{low matrix}\\
S^{up}  &= X^{up} + \sigma^{up}I\label{up matrix}
\end{align}
where column spaces of matrices $X^{lo}$ and $X^{up}$ are spanned by $\boldsymbol {x^{lo}}$ and $\boldsymbol {x^{up}}$ respectively. $\sigma^{lo}$ and $\sigma^{up}$ are scalars which are multiplied with identity matrix $I$ and added with column vectors of $X^{lo}$ and $X^{up}$ so as to form $S^{lo}$ and $S^{up}$. 
This structure of matrices, $S^{lo}$ and $S^{up}$ provides an alternative way to define constraint simplices, $\Delta_k^{lo}$ and $\Delta_k^{up}$ as,
\begin{align}
\Delta_k^{lo} &= \{\boldsymbol x = S^{lo}\boldsymbol\eta, \text{ }\boldsymbol\eta \in \mathds{R}^k|\sum_{i = 1}^{k}\eta_i = 1, \eta_i \geq 0\}\label{low alternate simplex}\\
\Delta_k^{up} &= \{\boldsymbol x = S^{up}\boldsymbol\xi, \text{ }\boldsymbol\xi \in \mathds{R}^k|\sum_{i = 1}^{k}\xi_i = 1, \xi_i \geq 0\}\label{up alternate simplex}
\end{align}

Hence it is possible to represent known population state $\boldsymbol x$ onto $\Delta_k^{lo}$ and $\Delta_k^{up}$ in terms of $\boldsymbol{\eta}$ and $\boldsymbol{\xi}$ respectively. 
\begin{align}
\boldsymbol{\eta} &= (S^{lo})^{-1} \boldsymbol x = \frac{1}{\sigma^{lo}}\{I - X^{lo}\}\label{ETA}\\
\boldsymbol{\xi} &= (S^{up})^{-1} \boldsymbol x = \frac{1}{\sigma^{up}}\{I - X^{up}\}\label{PHI}
\end{align}

As $\sum_{i = 1}^{k} x_i = 1$, \eqref{ETA} and \eqref{PHI} can be interpreted in simplified form as,
\begin{align}
\eta_i = \frac{1}{\sigma^{lo}}(x_i - x_i^{lo})\label{eta}\\
\xi_i = \frac{1}{\sigma^{up}}(x_i - x_i^{up})\label{phi}
\end{align} 

Here, $\boldsymbol{\eta}$ which correlates with lower bounds $\boldsymbol{x^{lo}}$ is increasing in nature while $\boldsymbol{\xi}$ is monotonically decreasing in nature which corresponds to upper bounds $\boldsymbol{x^{up}}$. It is so because $\sigma^{lo} > 0$ and $\sigma^{up} < 0$, which ensures that intersection of $\Delta_k^{lo}$ and $\Delta_k^{up}$ is not empty, \cite{ovalle2017escort}.

To incorporate the both lower and upper bounds \eqref{local_simplex_bound} the final escort function formulation consists multiplication of \eqref{eta} and \eqref{phi}. 
\begin{equation}
    \phi(x_i) = \eta_i \times \xi_i, \text{ for all } i = 1, 2, ..., k \label{final_escort_function}
\end{equation}

The escort function formulation \eqref{final_escort_function} ensures that if the initial distribution estimate belongs to the interior of the constraints \eqref{local_simplex_bound} then the individual strategy proportion will always remain bounded during the evolution process.

\section{Escort Dynamical approach for BTC problem}

The BTC temperature dynamics \eqref{temp_dynamics} with objective function \eqref{local_objective} depicts the plant dynamics \eqref{error_plant_module} whereas the DED dynamics \eqref{DED} represents the controller model $q(\cdot)$ in \eqref{controller}. 

The \textit{Assumption 4} demands controller dynamics \eqref{controller} to subdue as system state dynamics \eqref{plant} attains equilibrium point. If $\boldsymbol{e_t} = 0$, i.e. $\boldsymbol{t} - \boldsymbol{t^*} = 0$, then according to \eqref{local_objective} the respective output vector $\boldsymbol{f} = 0$. This in turn will diminish the controller dynamics given by \eqref{DED}, i.e. $e_{x_i} = 0$ for all $i = 1, 2, ..., k$. Hence, BTC dynamics \eqref{temp_dynamics} controlled by proposed DED protocol \eqref{DED} satisfy \textit{Assumption 4}.

The connectivity \textit{Assumption 3} of the system inherent communication topology ensures the fulfillment of C1 in \textit{Theorem 1}. The dynamics \eqref{temp_dynamics} is shown to be satisfying condition C2 in \cite{obando2015distributed}. Moreover, according to the positive invariantness property \eqref{positive invariantness} the DED dynamics \eqref{DED} satisfies C3. The fulfillment of C3 ensures that the DED model is passive in nature. Hence, the interconnection of two passive dynamics \eqref{temp_dynamics} and \eqref{DED} results in the stable rest point $(t^*, x^*)$ which corresponds to the attainment of output consensus.

Hence, according to \textit{theorem 1} the BTC temperature dynamics \eqref{temp_dynamics} for individual room with respective objective \eqref{local_objective}, driven by consensus like DED dynamic \eqref{DED} will reach equilibrium point obtaining output consensus i.e. welfare scheme ensuring equal payoff to every room.

The BTC process is a dynamic in nature where at every instance the temperature values gets affected by thermal disturbance $\boldsymbol{d}$ and surrounding temperature $\boldsymbol{t^a}$ according to \eqref{temp_dynamics}. This temperature is then regulated around setpoint temperature $\boldsymbol{t^{set}}$ using actuator mechanism which is driven by DED dynamics \eqref{DED}.The algorithmic implementation of BTC problem using DED dynamics is provided in Algorithm 1.
\rule{0.485\textwidth}{0.5pt}
\textbf{1. BTC using DED approach}\\
\rule{0.485\textwidth}{0.5pt}
\begin{algorithmic}[1]

\STATE \textbf{Initialize:} { $\boldsymbol{t^{set}}, \boldsymbol{t^{a}}, \boldsymbol{t^1}, \boldsymbol{W_z}, \boldsymbol{W_a}, \boldsymbol{X}, \boldsymbol{x^1},$ \\$\boldsymbol{x^{lo}}, \boldsymbol{x^{up}}, x^{total}$}

\STATE \textbf{Define:} {$\sigma^{lo} = x^{total} - \sum_{j = 1}^k{x^{lo}_j}$}\\
\hspace*{1.2cm}{$\sigma^{up} = x^{total} - \sum_{j = 1}^k{x^{up}_{j}}$}

\STATE \textbf{Compute:} $\eta_j = \frac{1}{\sigma^{lo}}(x_j - x^{lo}_j),$\\
\hspace*{1.2cm}$\xi_j = \frac{1}{\sigma^{up}}(x_j - x^{up}_j),$\\
\hspace*{1.2cm}$\phi(x_j) = \eta_j\times \xi_j$ for all $j = 1, 2, ..., k$ 

\STATE \textbf{Compute:} {$\boldsymbol{f^i} = [f_1^i, f_2^i, ..., f_k^i]^T$}\\
\hspace*{1.2cm}where, $f_j^i = t_j^i - t_j^{set_i}$

\STATE \textbf{Evaluate:} {$\boldsymbol{\dot{x}^{i}} = [\dot{x}_1^{i}, \dot{x}_2^{i}, ..., \dot{x}_k^{i}]$}\\
\hspace*{1.2cm}where, $\dot{x}_j^{i}$ is evaluated according to \eqref{DED}

\STATE \textbf{Compute:} {$\boldsymbol{x^{i+1}} = [x_1^{i+1}, x_2^{i+1}, ..., x_k^{i+1}]$}\\
\hspace*{1.2cm}where, $x_j^{i+1} = x_j^{i} + \dot{x}_j^{i}$

\STATE \textbf{Evaluate:} {$\boldsymbol{\dot{t}^i} = [\dot{t}_1^i, \dot{t}_2^i, ..., \dot{t}_N^i]$}\\
\hspace*{1.2cm}where, $\dot{t}_j^i$ is evaluated according to \eqref{temp_dynamics}

\STATE \textbf{Compute:} {$\boldsymbol{t^{i+1}} = [t_1^{i+1}, t_2^{i+1}, ..., t_N^{i+1}]$}\\
\hspace*{1.2cm}where, $t_j^{i+1} = t_j^{i} + \dot{t}_j^i$

\STATE \textbf{Assign:} $i\leftarrow {i+1}$
\label{main_algo}
\end{algorithmic}
\rule{0.485\textwidth}{0.5pt}

\section{REPRESENTATIVE CASE STUDY AND RESULTS}
The BTC issue formulated in Section II has been addressed using the DED approach. For comparison purpose the same problem has been implemented using consensus-based resource allocation protocol of DIP, refer Appendix-A for its mathematical formulation.

\subsection{Operational Scenario}
\begin{figure}[ht]
\centering
\begin{subfigure}{.24\textwidth}
  \centering
  \includegraphics[width=1\linewidth]{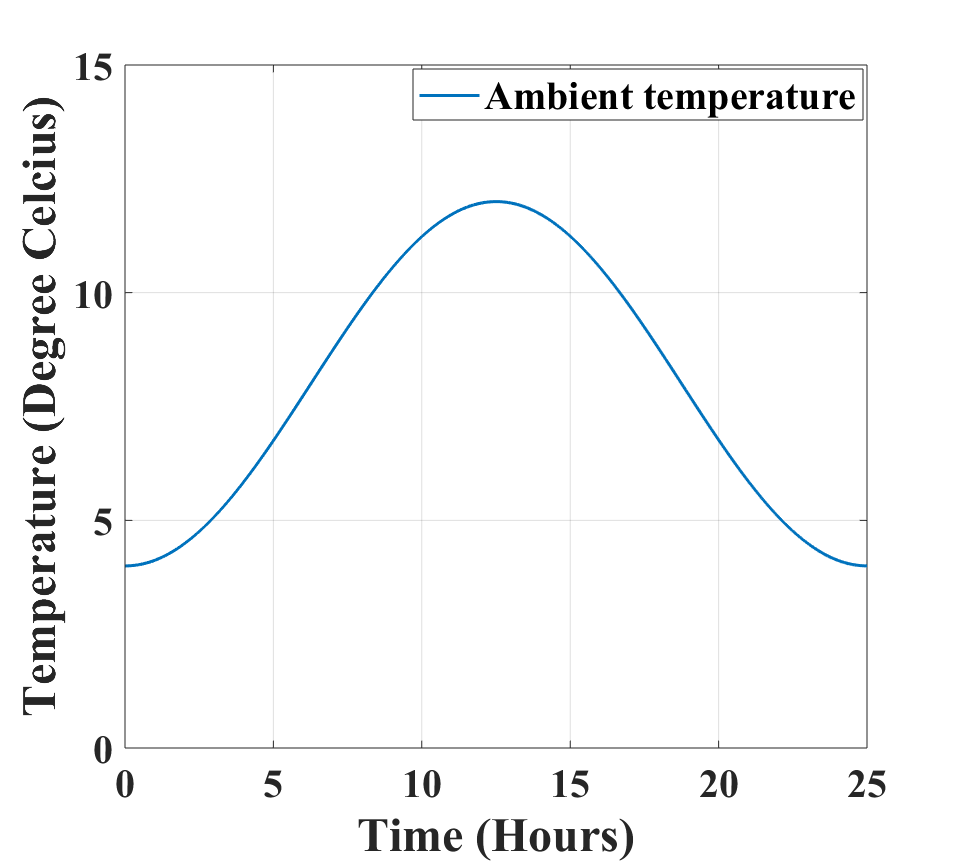}
  \caption{}
  \label{Ambient Temperature}
\end{subfigure}%
\begin{subfigure}{.24\textwidth}
  \centering
  \includegraphics[width=1\linewidth]{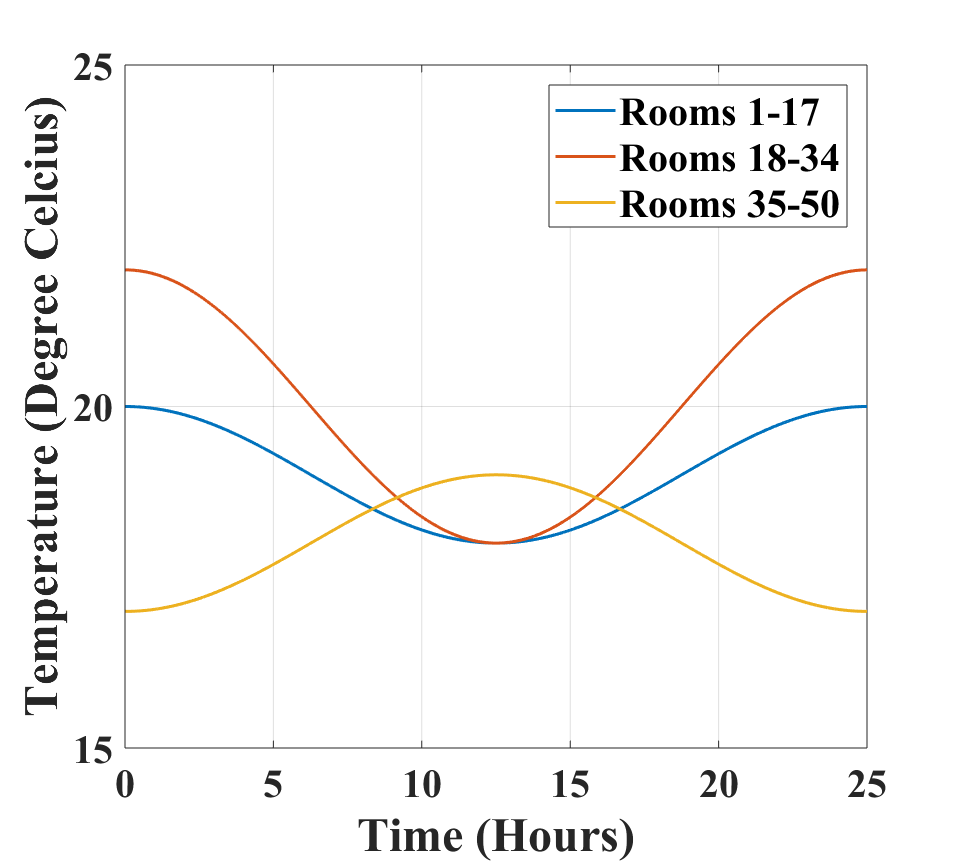}
    \caption{}
  \label{Desired Temperature}
\end{subfigure}
\caption{Operational conditions (a) Ambient temperature profile (b) Desired room temperature profiles}
\end{figure}

For implementation purpose building with $50$ rooms surrounded by the comparatively cold environment is analyzed over a period of an entire day, i.e. $k = 50$. The ambient temperature of the surrounding is varying as shown in Fig. \ref{Ambient Temperature}. Fig. \ref{Desired Temperature} depicts the desired temperature trajectories for the room temperatures throughout the day, where $50$ rooms are sectioned into $3$ groups viz., $1-17, 18-34, 35-50$. It is possible to define a separate temperature profile for each individual room however to neatly display the performance the rooms are divided into $3$ groups. Initially, every room is at $13^0 C$ temperature whereas initial actuator output is considered to be $0.5 \text{ kWh}$. Every actuator output is locally constrained within the range of $0 - 3.25 \text{ kWh}$ while global constrained is obtained by restricting the cumulative output of all actuators $130 \text{ kWh}$ at every given instance.

\subsection{Observations}
\begin{figure}[ht]
\centering
\includegraphics[width=0.41\textwidth]{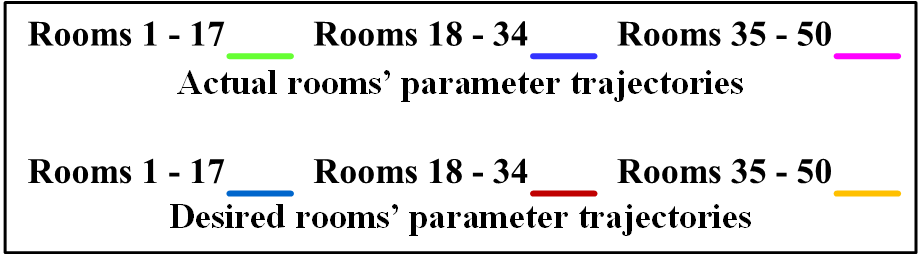}
\caption{Colour notations of actual and desired room parameters' trajectories}
\label{Legends}
\end{figure}

Observations include the comparison between desired trajectories and actually obtained trajectories of 50 rooms sectioned in 3 groups. The colour notations used to represent these trajectories is shown in Fig. \ref{Legends}.

\begin{figure}[h]
\centering
\begin{subfigure}{.235\textwidth}
  \centering
  \includegraphics[width=1\linewidth]{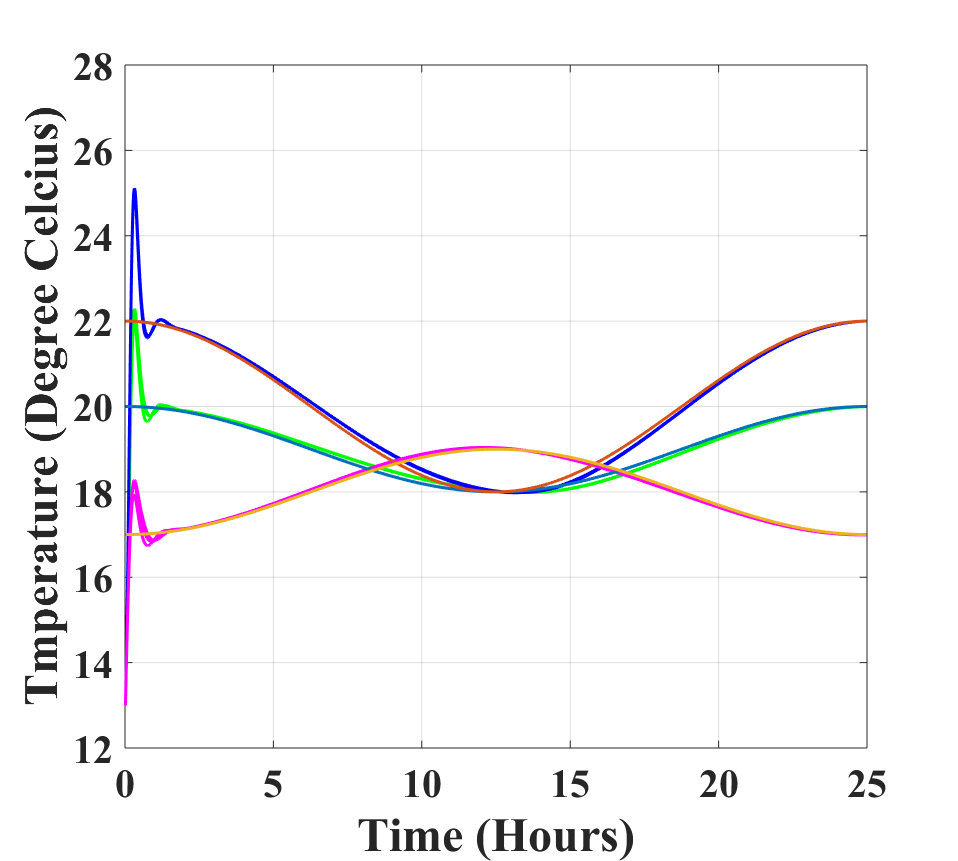}
  \caption{}
  \label{DED Temperature}
\end{subfigure}%
\begin{subfigure}{.235\textwidth}
  \centering
  \includegraphics[width=1\linewidth]{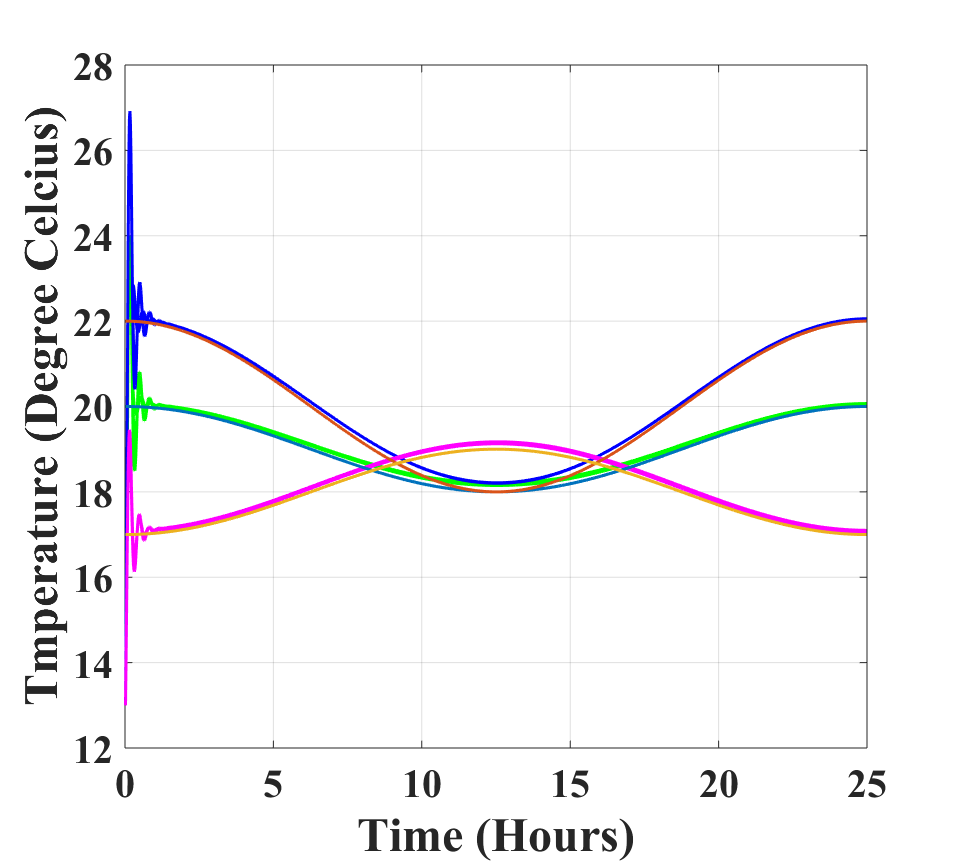}
  \caption{}
  \label{DIP Temperature}
\end{subfigure}
\caption{Desired temperature trajectories plotted against temperature trajectories obtained by implementing (a) DED approach (b) DIP approach}
\end{figure}

\begin{figure}[ht]
\centering
\begin{subfigure}{.235\textwidth}
  \centering
  \includegraphics[width=1\linewidth]{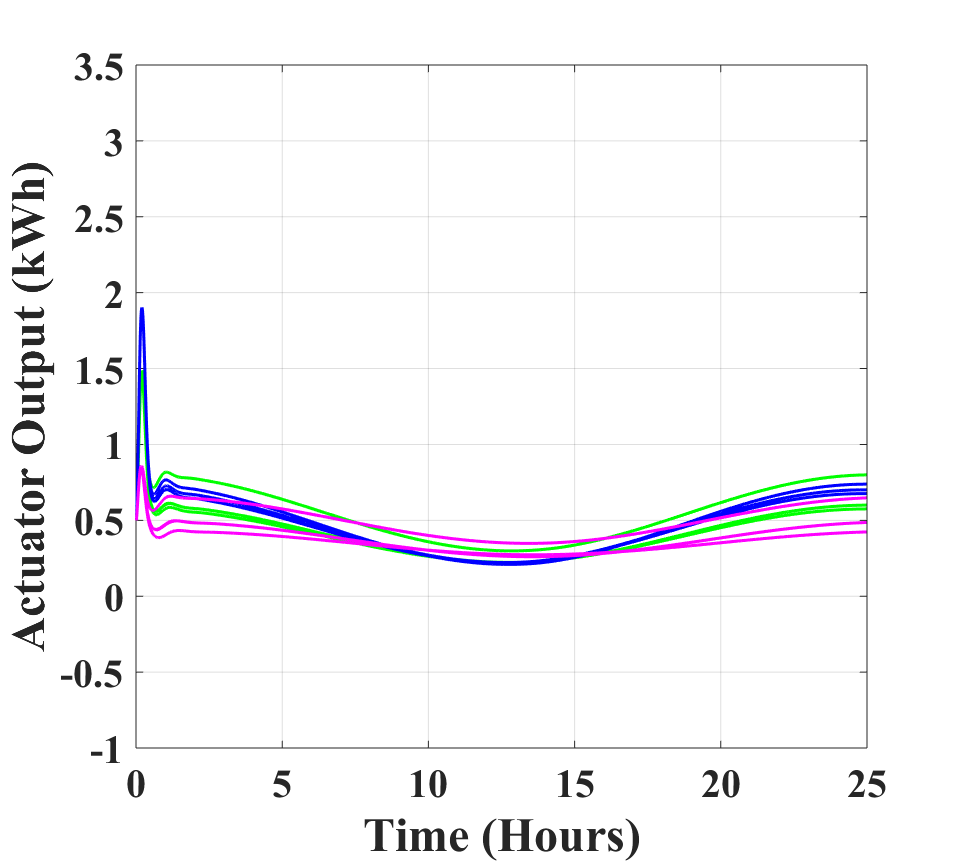}
  \caption{}  
  \label{DED Actuator}
\end{subfigure}%
\begin{subfigure}{.235\textwidth}
  \centering
  \includegraphics[width=1\linewidth]{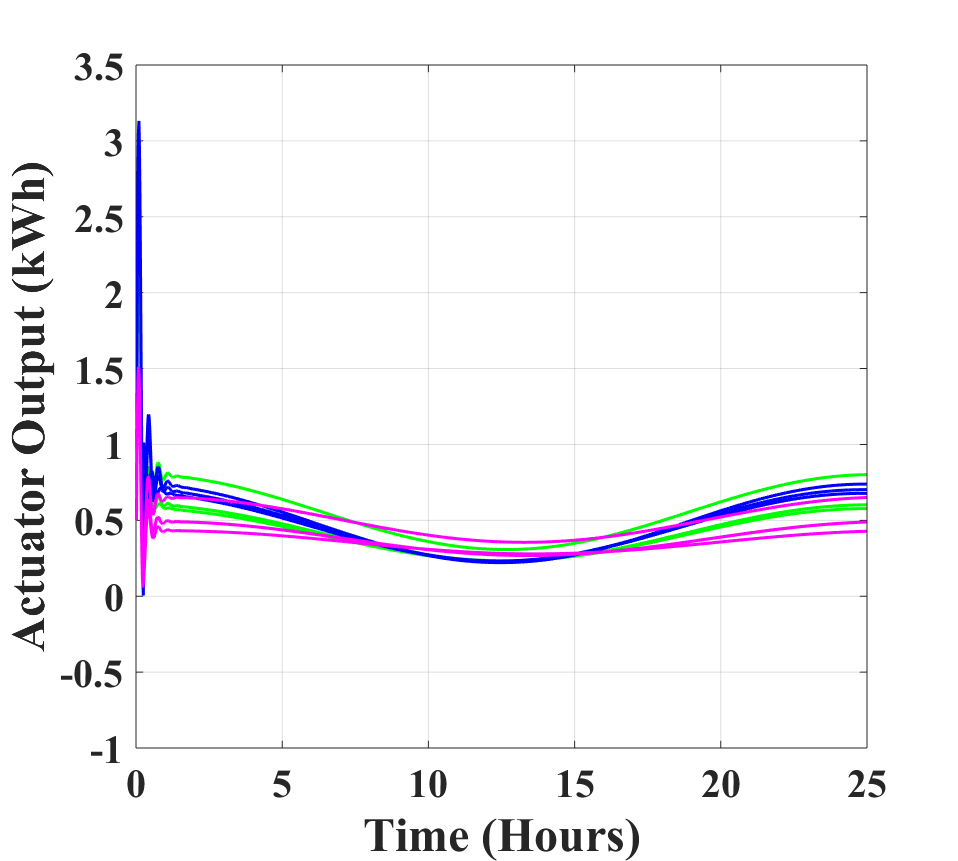}
  \caption{}
  \label{DIP Actuator}
\end{subfigure}
\caption{Actuator trajectories obtained by implementing (a) DED approach (b) DIP approach}
\label{Actuator output}
\end{figure}

\begin{figure}[h]
\centering
\begin{subfigure}{.24\textwidth}
  \centering
  \includegraphics[width=1\linewidth]{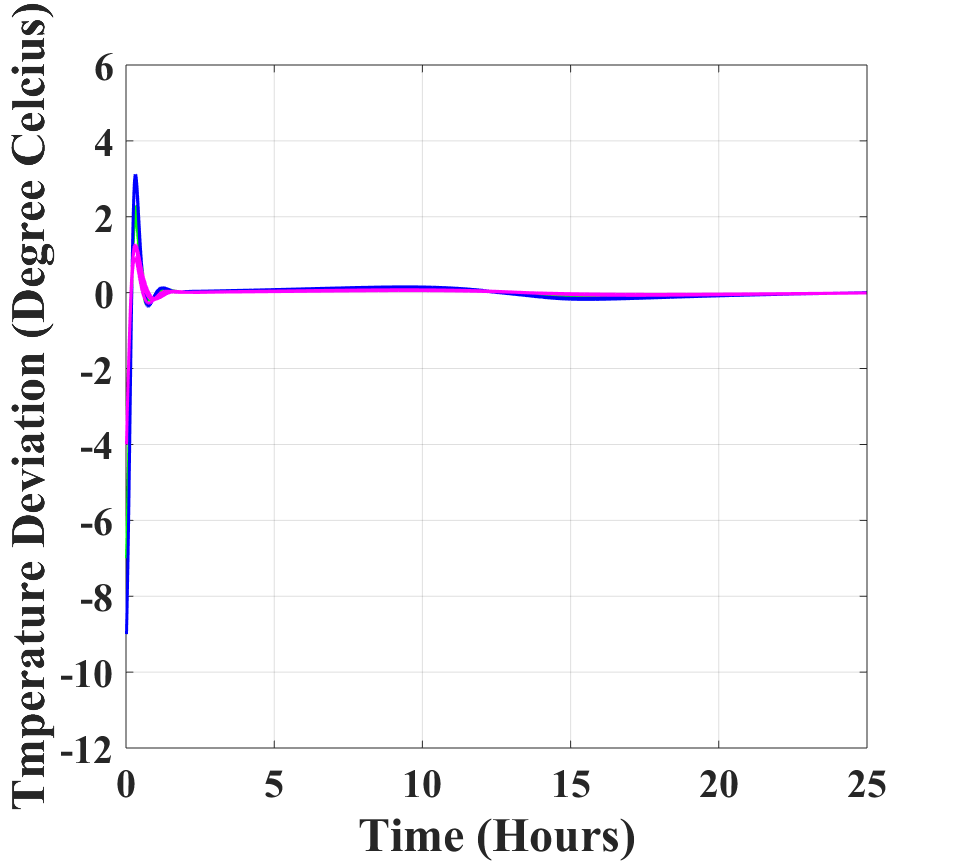}
  \caption{}
  \label{DED Payoff}
\end{subfigure}%
\begin{subfigure}{.24\textwidth}
  \centering
  \includegraphics[width=1\linewidth]{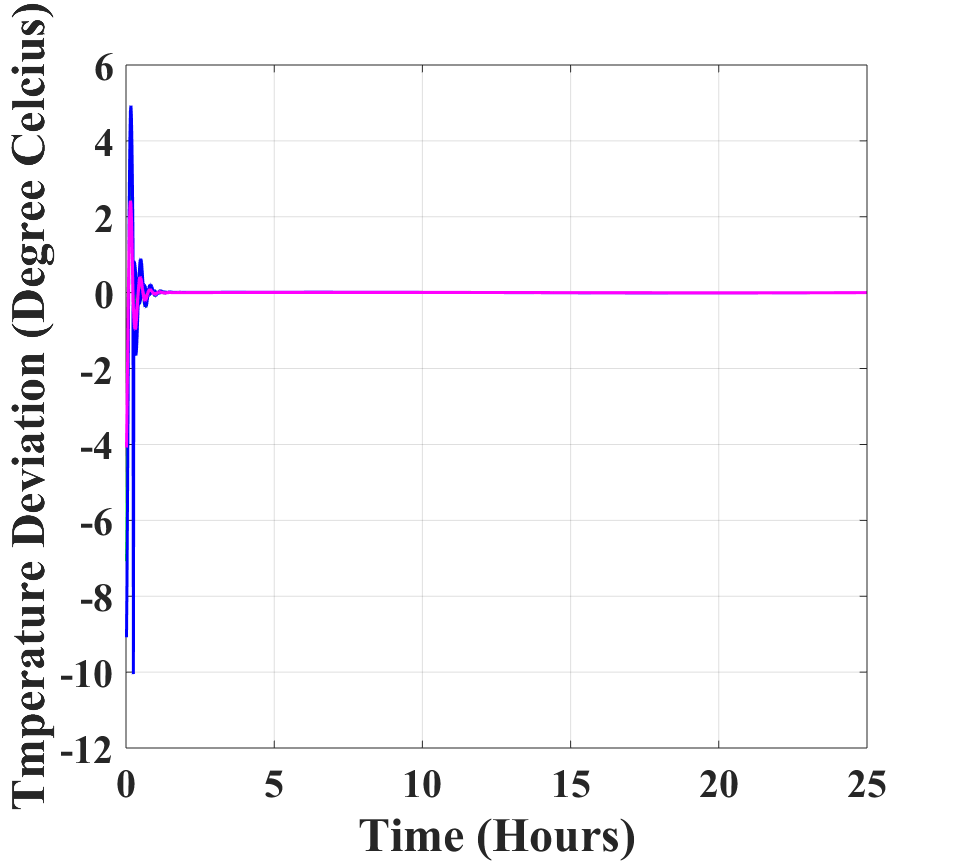}
  \caption{}
  \label{DIP Payoff}
\end{subfigure}
\caption{Payoff trajectories computed for (a) DED approach (b) DIP approach}
\end{figure}

Results obtained for implementation of building temperature control using DED are represented in Fig. \ref{DED Temperature},  \ref{DED Actuator}, and \ref{DED Payoff}. In the initial operational phase fast increase in actuator outputs is observed which results in on average $2^0 C$ temperature increment above desired set points. Corresponding oscillations can also be observed in payoff values when payoff values converge to zero, temperature profiles at respective interval attain desired temperature values, Fig. \ref{DED Payoff}. Within one damped oscillatory cycle temperature deviations are reduced and desired trajectories are traced thereafter Fig. \ref{DED Temperature}. After initial perturbations, actuator profiles display smooth variations, Fig. \ref{DED Actuator}.

Fig. \ref{DIP Temperature},  \ref{DIP Actuator}, and \ref{DIP Payoff} represent observations regarding implementation of the DIP method for the same problem. Like DED, DIP displayed oscillatory behaviour in the initial phase, also referred to as startup transience, but it persisted for more than one oscillation cycles, Fig. \ref{DIP Actuator}. Similar fluctuations are reflected in initial temperature profiles, as temperature varies back and forth of the desired value, Fig. \ref{DIP Temperature}. Peak temperature, as well as actuator output overshoots in case of DIP, are higher than that of DED, which results in larger temperature deviations shown in Fig. \ref{DIP Payoff}.

Also, temperature trajectories utilizing DED approach are closer to desired temperature trajectories than that of trajectories implementing DIP. When ambient temperature of the surrounding is above the average value the corresponding actuator outputs are low, Fig. \ref{Actuator output}. It has been observed that the initial peak overshoot is directly proportional to the initial difference between the desired and actual temperature.

\subsection{Analysis}
In an ideal implementation algorithm, real-time trajectories should convergence fast to desired temperature trajectories without any overshoots. DED, DIP being a step size dependent first-order dynamical protocols show some initial oscillatory behaviour.
 
\subsubsection{Presence of startup transience} 
\begin{figure}[ht]
\centering
\begin{subfigure}{.25\textwidth}
  \centering
  \includegraphics[width=1\linewidth]{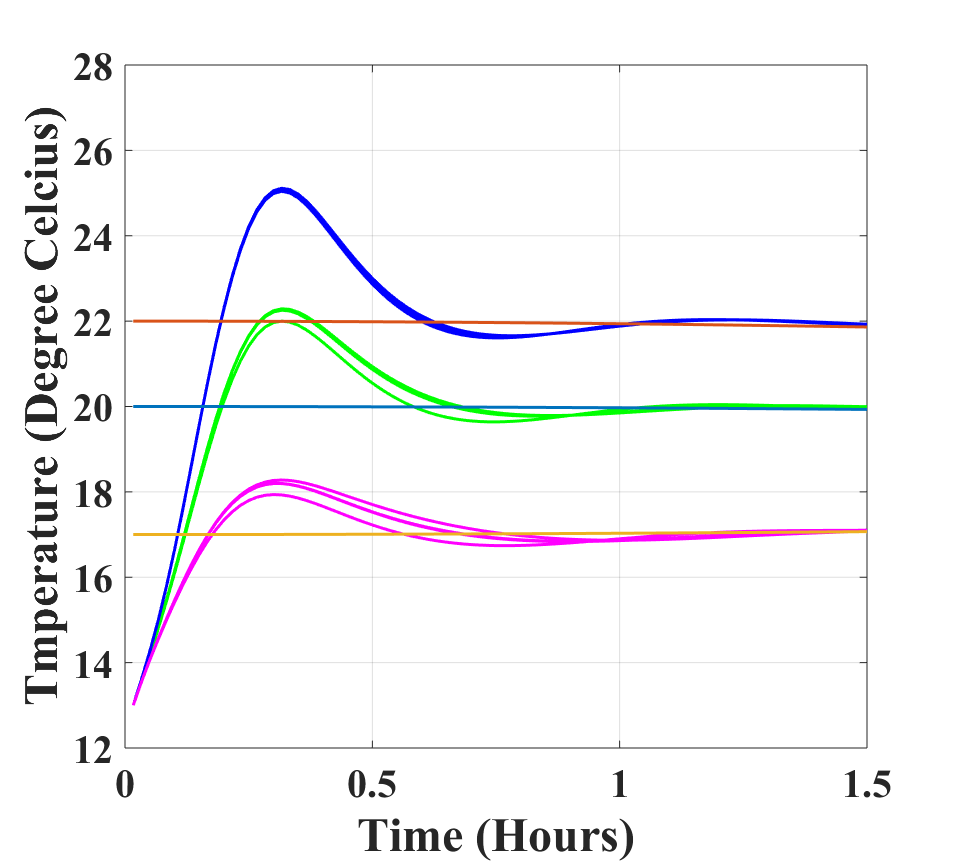}
  \caption{}
  \label{Zoom DED Temperature}
\end{subfigure}%
\begin{subfigure}{.25\textwidth}
  \centering
  \includegraphics[width=1\linewidth]{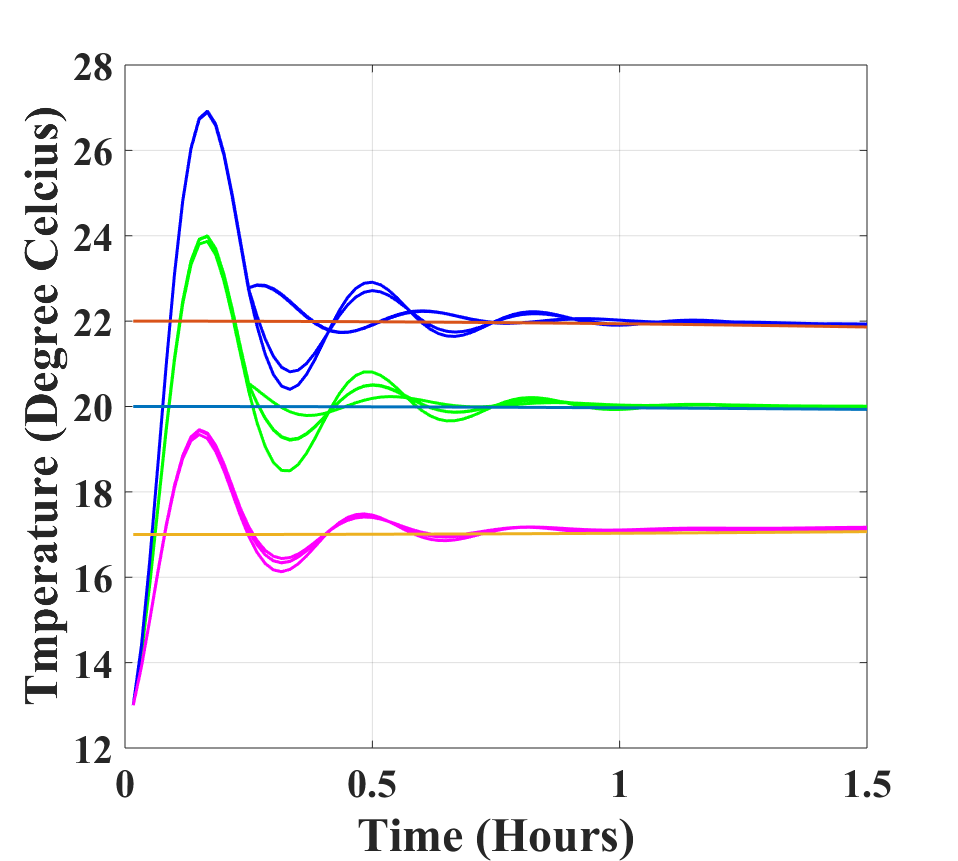}
  \caption{}
  \label{Zoom DIP Temperature}
\end{subfigure}
\caption{Desired temperature trajectories plotted against temperature trajectories obtained by implementing (a) DED approach (b) DIP approach}
\label{Temperature comparison}
\end{figure}

\begin{figure}[ht]
\centering
\begin{subfigure}{.25\textwidth}
  \centering
  \includegraphics[width=1\linewidth]{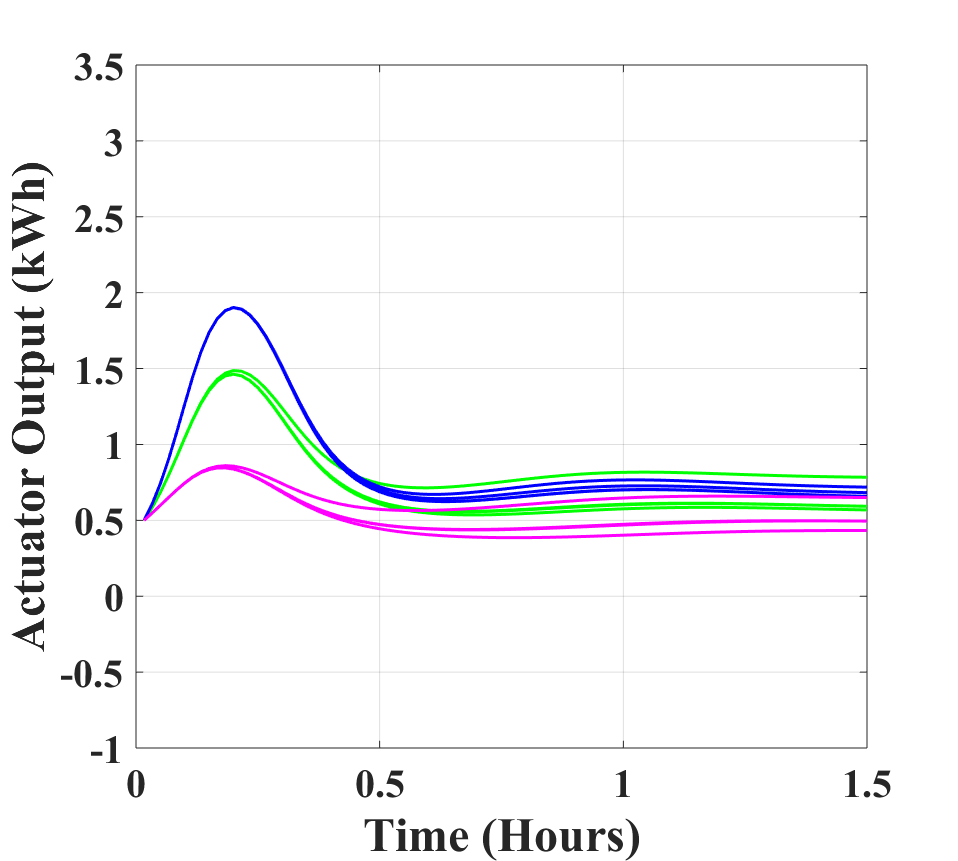}
  \caption{}  
  \label{Zoom DED Actuator}
\end{subfigure}%
\begin{subfigure}{.25\textwidth}
  \centering
  \includegraphics[width=1\linewidth]{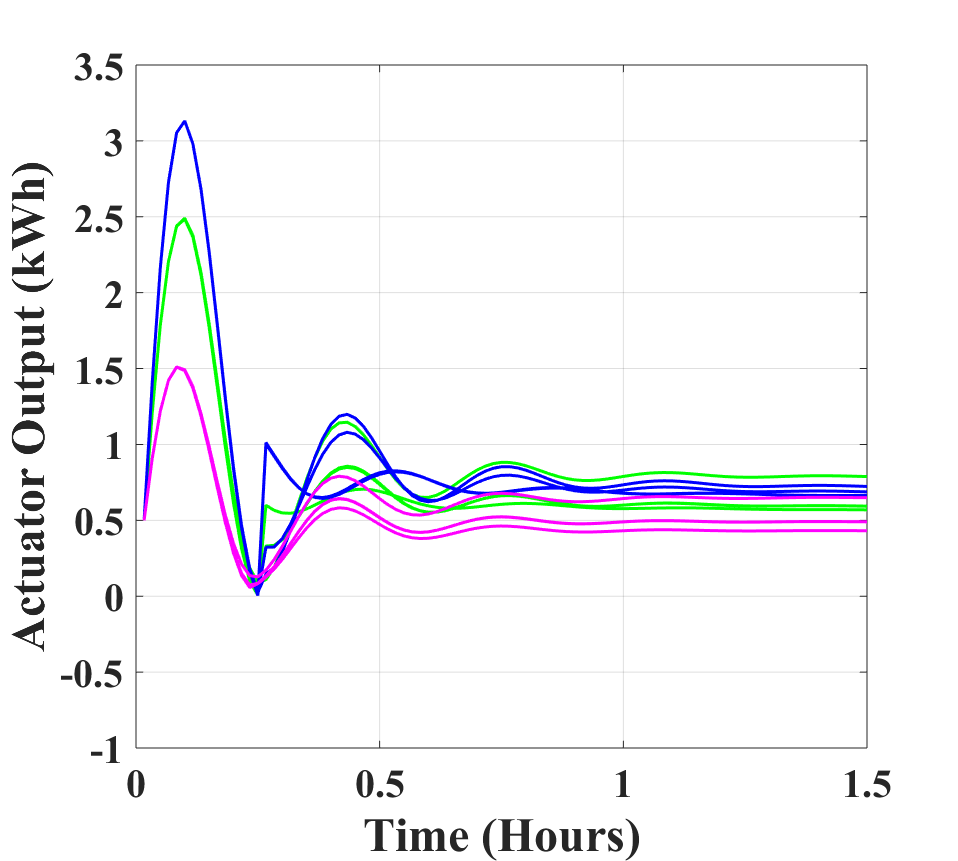}
  \caption{}
  \label{Zoom DIP Actuator}
\end{subfigure}
\caption{Actuator trajectories obtained by implementing (a) DED approach (b) DIP appraoch}
\label{Actuator comparison}
\end{figure}

The oscillatory response for DED approach is less as compared to DIP approach, Fig. \ref{Temperature comparison} and \ref{Actuator comparison}. In DIP approach, control parameters i.e. actuator outputs are constrained by means of a barrier function, \eqref{barrier function}, which actively participates in manipulating actuator parameter, keeping it in the centre of the predefined range. When actuator value reaches toward one of the boundaries the payoff value changes drastically which in turn force it towards the centre. This operation introduces oscillatory behaviour over a couple of cycles, Fig. \ref{Zoom DIP Temperature} and \ref{Zoom DIP Actuator}. On the other hand, the inclusion of constraint directly in the dynamical equation \eqref{ED} using escort function guarantees a comparatively smooth transition from initial conditions to the desired ones, without violating the constraints, Fig. \ref{Zoom DED Actuator} and \ref{Zoom DED Actuator}. This is because payoff values are not directly affected by changes in actuator output. 

\subsubsection{Desired trajectory tracking} 
In DED approach, deviations in temperature value is defined as a payoff. As operations continue every sub-system (room) reach on a common consensus value of payoff function. The negative difference between actual and required temperature value increases the actuator output value and vice versa. Hence algorithm implementing DED tracks desired trajectory closely, Fig. \ref{DED Temperature}. While in DIP approach along with temperature deviations cost function also includes barrier function. Hence when DIP reaches consensus, the consensus value gets biased by the barrier function value, which introduces a marginal error between actual temperature trajectories and desired trajectories, Fig. \ref{DIP Temperature}.

Lower actuator output value in comparatively warm surrounding environment, Fig. \ref{Actuator output} corresponds to efficient energy utilization. 

\section{Conclusion}
The issue of BTC in the resource allocation framework has been analyzed using evolutionary game theory. The distributed version of evolutionary dynamics of ED has been utilized to address resource allocation problem within system-defined constraints. The local bounds on the individual actuator bounds are incorporated within the DED framework through the concept of the intersection of simplices. The consensus like DED approach is shown to attain output consensus for dynamic resource allocation problem with global resource constraint. The performance analysis on the metric of smooth trajectory tracking and low startup transience of BTC mechanism implementation through DED approach and DIP approach is carried out. Unlike DED where local bounds on individual actuator output are accommodated through barrier formulation, the DED approach restricts the dynamics within the dynamics through escort function formulation. This escort function formulation is shown to have better performance as compared to DIP approach in terms of mechanism longevity through low startup transience than that of DIP.

This research work can be extended to dynamic resource allocation problems where the constraint set is time-varying. This can be incorporated through varying domains of the intersection of simplices.

\section*{Appendix}
\subsection{Distributed Interior Point method}


Considering the BTC problem described in Section II, DIP protocol used to obtain next control vector values, is given as,
\begin{equation}
\dot v_i = \sum_{j\in \mathcal{M}_i}[r_j(\boldsymbol v) - r_i(\boldsymbol v)]; \label{DIP}
\end{equation}
where, $(1 \times k)$ dimensional payoff vector is constructed as
\begin{equation}
r_i(\boldsymbol v) = f_i(\boldsymbol v) + b_i(v_i) \text{ for } i = 1, 2, ..., k\label{cost function} 
\end{equation}
here, $b_i(\boldsymbol v_i)$ is defined as the derivative of barrier function $\mathcal{B}(\boldsymbol v)$ with respect to $v_i$ is added with original objective of the $i^{th}$ sub-system to form payoff function for respective sub-system.

Barrier function restricts the control vector values within a predefined constraints. Logarithmic barrier function is one of the barrier function, \cite{obando2015distributed}, given as
\begin{equation}
\mathcal{B}(\boldsymbol v) = -\epsilon[in(v_i - v_i^{lo}) + ln(v_i^{up} - v_i)], \forall i = 1, 2, ..., k\label{barrier function}
\end{equation}

\bibliographystyle{IEEEtran}
\bibliography{References}

\end{document}